\documentclass[11pt]{article}%
\usepackage{amsmath}
\usepackage{amsfonts}
\usepackage{amssymb}
\usepackage{graphicx}%
\providecommand{\U}[1]{\protect\rule{.1in}{.1in}}
\newtheorem{theorem}{Theorem}[section]

\newtheorem{remark}[theorem]{Remark}

\numberwithin{equation}{section}

\begin{document}

\title{Reconstruction techniques for inverse Sturm-Liouville problems with complex coefficients}
\author{Vladislav V. Kravchenko\\{\small Departamento de Matem\'{a}ticas, Cinvestav, Unidad Quer\'{e}taro, }\\{\small Libramiento Norponiente \#2000, Fracc. Real de Juriquilla,
Quer\'{e}taro, Qro., 76230 MEXICO.}\\{\small e-mail: vkravchenko@math.cinvestav.edu.mx}}
\maketitle

\begin{abstract}
A variety of inverse Sturm-Liouville problems is considered, including the
two-spectrum inverse problem, the problem of recovering the potential from the
Weyl function, as well as the recovery from the spectral function. In all
cases the potential in the Sturm-Liouville equation is assumed to be complex
valued. A unified approach for the approximate solution of the inverse
Sturm-Liouville problems is developed, based on Neumann series of Bessel
functions (NSBF) representations for solutions and their derivatives. Unlike
most existing approaches, it allows one to recover not only the complex-valued
potential but also the boundary conditions of the Sturm-Liouville problem.
Efficient accuracy control is implemented. The numerical method is direct. It
involves only solving linear systems of algebraic equations for the
coefficients of the NSBF representations, while eventually the knowledge only
of the first NSBF coefficients leads to the recovery of the Sturm-Liouville
problem. Numerical efficiency is illustrated by several test examples.

\end{abstract}

\section{Introduction}

Let $q\in%
\mathcal{L}%
^{2}\left(  0,b\right)  $ be a complex-valued function. Consider the
Sturm-Liouville equation%
\begin{equation}
-y^{\prime\prime}+q(x)y=\lambda y, \label{Schr Intro}%
\end{equation}
on a finite interval $0<x<b$.

In the present work, we discuss a method for approximately solving a range of
inverse spectral problems associated with equation (\ref{Schr Intro}). In
particular, (i) the two-spectrum inverse problem; (ii) the recovery of the
Sturm-Liouville problem from its Weyl function; (iii) the recovery of the
Sturm-Liouville problem from a spectrum and corresponding multiplier
constants, which represent the values of normalized eigenfunctions at
endpoints; (iv) the recovery of the Sturm-Liouville problem from a spectrum
and corresponding norming constants (weight numbers). While the method can be
extended to more general situations, for simplicity, throughout this paper, we
assume that the spectra considered are simple.

Development of efficient methods for solving inverse Sturm-Liouville problems
is an active research area. Several methods have been proposed for the
solution of inverse Sturm-Liouville problems of different types, which do not
involve the techniques presented below (see \cite{Andersson1990},
\cite{BDK2005}, \cite{Bockmann et al 2011}, \cite{Brown et al 2003},
\cite{Drignei 2015}, \cite{Gao et al 2013}, \cite{Gao et al 2014},
\cite{IgnatievYurko}, \cite{Kammanee Bockman 2009}, \cite{Lowe et al 1992},
\cite{Neamaty et al 2017}, \cite{Neamaty et al 2019}, \cite{Rohrl},
\cite{Rundell Sacks}, \cite{Sacks} and many other publications). However,
usually they require the knowledge of additional parameters like, e.g., the
parameter $\omega:=\frac{1}{2}\int_{0}^{b}q(t)\,dt$ (such is the case of
\cite{IgnatievYurko}, \cite{Rundell Sacks}, \cite{Sacks}) and do not allow the
recovery of the boundary conditions. As an illustration of this fact we quote
\cite[p. 177]{Rundell Sacks}: \emph{`as was mentioned earlier, if complete
spectral data is available, then it is in theory possible to determine the
boundary conditions as part of the solution of the problem. We do not believe,
however, that this is numerically feasible in most cases.'} For this reason
the boundary conditions are supposed to be known in all the above mentioned
publications. Our approach is free of these restrictions, and moreover, to the
difference from those previous publications, it allows us to recover complex
valued potentials and complex constants from boundary conditions.

The overall approach presented here is based on Neumann series of Bessel
functions (NSBF) representations for solutions of (\ref{Schr Intro}) and for
their derivatives. The NSBF representations were derived in \cite{KNT}, where
results on their convergence and on the decay rate for their coefficients can
be found. In several publications, these NSBF representations have been
applied successfully to solving inverse spectral problems. First, in
\cite{Kr2019JIIP}, the classical problem of recovering the Sturm-Liouville
problem from its spectral function (see problem (IP4) below) was approached
with the aid of the NSBF and the Gelfand-Levitan integral equation. Due to the
discontinuity of the sum of the series representing the Gelfand-Levitan
integral kernel, the method provided less accuracy near endpoints. In
\cite[Sect. 13.4]{KrBook2020} its simple modification was proposed, solving
this particular problem. In \cite{KKT2019AMC} another series representation
for the Gelfand-Levitan kernel was derived, such that the sum of the series
became continuous. This was used in \cite{KT2021 IP1} to improve the approach
from \cite{Kr2019JIIP}. The approach based on using the NSBF combined
with\ the Gelfand-Levitan integral equation was applied to other inverse
Sturm-Liouville problems, both on a finite interval and on a half-line
\cite{DKK2019MMAS}, \cite{KKK QuantumFest}, \cite{KarapetyantsKravchenkoBook},
\cite[Sect. 13.4]{KrBook2020}, \cite{KST2020IP}, \cite{KT2021 IP1},
\cite{KT2021 IP2}, \cite{KV2022JMS2}.

Another way of applying NSBF representations to coefficient inverse problems
was explored in \cite{KKC2022Mathematics}, \cite{KrSpectrumCompl},
\cite{AvdoninKravchenko2023JIIP}, \cite{AKK2024IPI}, \cite{AKK2024MMAS},
\cite{CKK2024MMAS}. From the input data of the problem and by computing NSBF
coefficients (evaluated at $x=b$), a couple of the characteristic functions of
certain Sturm-Liouville problems were computed. They were used for computing a
spectrum and a sequence of corresponding multiplier constants (which relate
the linearly dependent eigenfunctions normalized at opposite endpoints), thus
reducing the original problem to problem (IP3) below, which was solved by
substituting NSBF into the relation between the linearly dependent
eigenfunctions and solving the resulting system of linear algebraic equations.
This scheme works well. The spectrum and the multiplier constants are computed
in great quantities and with a remarkable and uniform accuracy. It is also
interesting to note that this approach allows one to compute spectra of
Sturm-Liouville problems from other input data and without knowing the
potential $q(x)$. Moreover, one even can complete a spectrum, i.e., by given
several first eigenvalues of a Sturm-Liouville problem to compute much more of
them, again without knowing $q(x)$ \cite{KrSpectrumCompl}. Nevertheless, this
approach is restricted to real valued potentials, since computing complex
eigenvalues is a challenging task.

Another approach, which does not require computing eigenvalues, and therefore,
it is suitable for complex-valued potentials, was proposed in \cite{Kr2024JMP}
for solving a coefficient inverse problem for (\ref{Schr Intro}). For
obtaining the main system of linear algebraic equations for the NSBF
coefficients, instead of the relation between the linearly dependent
eigenfunctions another identity, relating three different solutions of
(\ref{Schr Intro}) was used (see identity (\ref{T=}) below). In the present
work we develop further this idea and devise an approach, which allows one to
solve a wide variety of inverse Sturm-Liouville problems for equation
(\ref{Schr Intro}) with a complex-valued potential $q(x)$, including the
recovery of the unknown constants in the boundary conditions. It is based on
the NSBF representations and several identities for solutions of
(\ref{Schr Intro}). The overall scheme reduces the solution of an inverse
problem to solving a couple of systems of linear algebraic equations. The
potential and constants are recovered from the very first NSBF coefficients.
The resulting numerical method is simple, accurate and fast. Several practical
criteria are derived for choosing the involved parameters and ensuring
accuracy. The method requires neither additional assumptions on the potential
nor initial guesses and does not involve any iterative procedures. Thus, the
contribution of this work consists of a universal approach to the approximate
solution of a wide range of inverse Sturm-Liouville problems for equation
(\ref{Schr Intro}).

\section{Two-spectrum inverse problem setting}

Let $q\in%
\mathcal{L}%
^{2}\left(  0,b\right)  $ be complex valued, $b>0$. Consider the
Sturm-Liouville equation%
\begin{equation}
-y^{\prime\prime}+q(x)y=\lambda y,\quad x\in(0,b), \label{Schr}%
\end{equation}
where $\lambda\in\mathbb{C}$ is a spectral parameter.

Denote by $L=L(q(x),h,H)$ the Sturm-Liouville problem for (\ref{Schr}) with
the boundary conditions
\begin{equation}
U(y):=y^{\prime}(0)-hy(0)=0,\quad V(y):=y^{\prime}(b)+Hy(b)=0, \label{bc1}%
\end{equation}
where $h$ and $H$ are some complex constants.

Denote by $L^{0}=L^{0}(q(x),H)$ the Sturm-Liouville problem for (\ref{Schr})
with the boundary conditions%
\begin{equation}
y(0)=0,\quad V(y)=0. \label{bc2}%
\end{equation}
Theory of problems $L$ and $L^{0}$ is well developed (see, \cite{Marchenko},
\cite{Brown et al 2002}, \cite{Buterin 2007}). In particular, the spectrum of
each problem is discrete, and a finite number of the eigenvalues can have
multiplicities greater than one. However, for simplicity, throughout this
paper we assume that the spectra considered are simple.

We start by considering the classical two-spectrum inverse problem, which is
formulated as follows (the other inverse Sturm-Liouville problems considered
in the present work, are introduced in Section
\ref{Sect Other inverse problems}).

\textbf{Inverse Problem (IP1) }Given two spectra $\left\{  \lambda
_{k}\right\}  _{k=0}^{\infty}$ and $\left\{  \lambda_{k}^{0}\right\}
_{k=0}^{\infty}$ of problems $L$ and $L^{0}$, respectively, find $q(x)\in%
\mathcal{L}%
^{2}\left(  0,b\right)  $ and the constants $h$, $H$.

For the uniqueness result for (IP1) we refer to \cite{Brown et al 2002},
\cite{Buterin 2007}, and for the stability result for (IP1) to
\cite{ButerinKuznetsova2019}.

The primary result of this study is a method for approximating the solution of
(IP1) when a finite number of eigenpairs $\left\{  \lambda_{k},\,\lambda
_{k}^{0}\right\}  _{k=0}^{K}$ are available. Often it is convenient to deal
with the square root of the spectral parameter: $\rho=\sqrt{\lambda}$. Since
it is not essential which of the square roots is chosen, we can always assume
that $\operatorname{Im}\rho\geq0$. Denote $\rho_{k}=\sqrt{\lambda_{k}}$ and
$\mu_{k}=\sqrt{\lambda_{k}^{0}}$ the square roots of the eigenvalues of
problems $L$ and $L^{0}$, respectively. Sometimes $\rho_{k}$ and $\mu_{k}$ are
called singular numbers of the respective Sturm-Liouville problems.

\section{Solutions and their series representations \label{Sect Prelim}}

Let $\varphi_{h}(\rho,x)$, $S(\rho,x)$, $\psi_{H}(\rho,x)$, $T(\rho,x)$ denote
the solutions of (\ref{Schr}) satisfying the respective initial conditions
\[
\varphi_{h}(\rho,0)=1,\quad\varphi_{h}^{\prime}(\rho,0)=h,\quad S(\rho
,0)=0,\quad S^{\prime}(\rho,0)=1,
\]%
\[
\psi_{H}(\rho,b)=1,\quad\psi_{H}^{\prime}(\rho,b)=-H,\quad T(\rho,b)=0,\quad
T^{\prime}(\rho,b)=1.
\]
For all $\rho\in\mathbb{C}$, we have $U(\varphi_{h})=0$ and $V(\psi_{H})=0$.
Denote
\[
\Delta\left(  \rho\right)  =W\left[  \psi_{H}(\rho,x),\,\varphi_{h}%
(\rho,x)\right]  ,
\]
where $W$ stands for the Wronskian (we recall that the Wronskian of any pair
of solutions of (\ref{Schr}) is independent of $x$). Substituting $x=0$ and
$x=b$, we find out that
\begin{equation}
\Delta\left(  \rho\right)  =V(\varphi_{h})=-U(\psi_{H}).
\label{equality charact}%
\end{equation}
The set of zeros of $\Delta\left(  \rho\right)  $ coincides with the set of
the singular numbers $\rho_{k}$ of problem $L$. Thus, $\Delta\left(
\rho\right)  $ is a characteristic function of problem $L$. Analogously, the
characteristic function of problem $L^{0}$ has the form%
\[
\Delta^{0}\left(  \rho\right)  =\psi_{H}(\rho,0)=S^{\prime}(\rho
,b)+HS(\rho,b).
\]
The following easily verifiable identities hold%
\begin{equation}
\psi_{H}(\rho,x)=\Delta^{0}\left(  \rho\right)  \varphi_{h}(\rho
,x)-\Delta\left(  \rho\right)  S(\rho,x) \label{psiH=}%
\end{equation}
and
\begin{equation}
T(\rho,x)=\varphi_{h}(\rho,b)S(\rho,x)-S(\rho,b)\varphi_{h}(\rho,x).
\label{T=}%
\end{equation}

We will use the following series representations for the solutions,$\ $derived
in \cite{KNT}.

\begin{theorem}
[\cite{KNT}]\label{Th NSBF} The solutions $\varphi_{h}(\rho,x)$, $S(\rho,x)$,
$\psi_{H}(\rho,x)$, $T(\rho,x)$ admit the following series representations
\begin{align}
\varphi_{h}(\rho,x)  &  =\cos\rho x+\sum_{n=0}^{\infty}(-1)^{n}g_{n}%
(x)\mathbf{j}_{2n}(\rho x),\label{phi NSBF}\\
S(\rho,x)  &  =\frac{\sin\rho x}{\rho}+\frac{1}{\rho}\sum_{n=0}^{\infty
}(-1)^{n}s_{n}(x)\mathbf{j}_{2n+1}(\rho x),\label{S NSBF}\\
\psi_{H}(\rho,x)  &  =\cos\left(  \rho\left(  x-b\right)  \right)  +\sum
_{n=0}^{\infty}(-1)^{n}\psi_{n}(x)\mathbf{j}_{2n}(\rho\left(  x-b\right)
),\label{psi NSBF}\\
T(\rho,x)  &  =\frac{\sin\left(  \rho\left(  x-b\right)  \right)  }{\rho
}+\frac{1}{\rho}\sum_{n=0}^{\infty}(-1)^{n}t_{n}(x)\mathbf{j}_{2n+1}\left(
\rho\left(  x-b\right)  \right)  . \label{T NSBF}%
\end{align}
where $\mathbf{j}_{k}(z)$ stands for the spherical Bessel function of order
$k$ ($\mathbf{j}_{k}(z):=\sqrt{\frac{\pi}{2z}}J_{k+\frac{1}{2}}(z)$, see,
e.g., \cite{AbramowitzStegunSpF}). For every $\rho\in\mathbb{C}$ the series
converge pointwise. For every $x\in\left[  0,L\right]  $ the series converge
uniformly in any strip of the complex plane of the variable $\rho$, parallel
to the real axis. In particular, the remainders of their partial sums
$\varphi_{h,N}(\rho,x):=\cos\rho x+\sum_{n=0}^{N}(-1)^{n}g_{n}(x)\mathbf{j}%
_{2n}(\rho x)$ and $S_{N}(\rho,x):=\frac{\sin\rho x}{\rho}+\frac{1}{\rho}%
\sum_{n=0}^{\infty}(-1)^{n}s_{n}(x)\mathbf{j}_{2n+1}(\rho x)$ admit the
estimates
\begin{equation}
\left\vert \varphi(\rho,x)-\varphi_{h,N}(\rho,x)\right\vert \leq
\frac{\varepsilon_{N}(x)\,\sinh(ax)}{a}\quad\text{and}\quad\left\vert \rho
S(\rho,x)-\rho S_{N}(\rho,x)\right\vert \leq\frac{\varepsilon_{N}%
(x)\,\sinh(ax)}{a} \label{estc2}%
\end{equation}
for any $\rho$ belonging to a strip $\left\vert \operatorname{Im}%
\rho\right\vert \leq a$, $a>0$, where $\varepsilon_{N}(x)$ is a positive
function tending to zero when $N\rightarrow\infty$. Analogous estimates are
valid for the remainders of the series (\ref{psi NSBF}) and (\ref{T NSBF}).

The first coefficients of the series have the form
\begin{equation}
g_{0}(x)=\varphi_{h}(0,x)-1,\quad s_{0}(x)=3\left(  \frac{S(0,x)}{x}-1\right)
, \label{g0s0}%
\end{equation}%
\begin{equation}
\psi_{0}(x)=\psi_{H}(0,x)-1,\quad t_{0}(x)=3\left(  \frac{T(0,x)}%
{x-b}-1\right)  , \label{psi0t0}%
\end{equation}
and the rest of the coefficients can be calculated following a simple
recurrent integration procedure.
\end{theorem}

Analogous series representations are also available for the derivatives of the
solutions \cite{KNT}. In particular, we will use the series representations%
\begin{equation}
\varphi_{h}^{\prime}(\rho,x)=-\rho\sin\left(  \rho x\right)  +\left(
h+\omega(x)\right)  \cos\left(  \rho x\right)  +\sum_{n=0}^{\infty}%
(-1)^{n}\gamma_{n}(x)\mathbf{j}_{2n}(\rho x), \label{phi prime NSBF}%
\end{equation}
and
\begin{equation}
S^{\prime}(\rho,x)=\cos\left(  \rho x\right)  +\frac{\omega(x)}{\rho}%
\sin\left(  \rho x\right)  +\frac{1}{\rho}\sum_{n=0}^{\infty}(-1)^{n}%
\sigma_{n}(x)\mathbf{j}_{2n+1}(\rho x), \label{S prime NSBF}%
\end{equation}
where
\[
\omega(x):=\frac{1}{2}\int_{0}^{x}q(t)\,dt,
\]
and the first coefficients have the form
\[
\gamma_{0}(x)=g_{0}^{\prime}(x)-\omega(x),\quad\sigma_{0}(x)=\frac{s_{0}%
(x)}{x}+s_{0}^{\prime}(x)-3\omega(x).
\]
Again, the rest of the coefficients can be calculated following a simple
recurrent integration procedure, and the series (\ref{phi prime NSBF}) and
(\ref{S prime NSBF}) converge uniformly in any strip of the complex plane of
the variable $\rho$, parallel to the real axis.

We recall that according to \cite[p. 522]{Watson}: \emph{`Any series of the
type}
\begin{equation}
\sum_{n=0}^{\infty}a_{n}J_{\nu+n}(z) \label{Neumann series}%
\end{equation}
\emph{is called a Neumann series, although in fact Neumann considered only the
special type of series for which }$\nu$\emph{ is an integer; the investigation
of the more general series is due to Gegenbauer.'}\ Since in functional
analysis, the term Neumann series is often used for refering to series of
another nature, the series (\ref{Neumann series}) are often referred to as
Neumann series of Bessel functions (NSBF) (see, e.g., \cite{Wilkins}). Thus,
we will refer to the series representations for the solutions and their
derivatives, which are used in the present work, as NSBF representations.

\begin{remark}
The NSBF representations are the result of the series expansions of the
transformation (transmutation) operator kernels into corresponding series in
terms of Legendre polynomials. Namely, as it is well known (see, e.g.,
\cite{Freiling and Yurko}, \cite{LevitanInverse}, \cite{Marchenko},
\cite{SitnikShishkina Elsevier}) there exist continuous functions
$\mathbf{G}_{h}(x,t)$ and $\mathbf{S}(x,t)$ in the domain $0\leq t\leq x\leq
b$, such that
\[
\varphi_{h}(\rho,x)=\cos\rho x+\int_{0}^{x}\mathbf{G}_{h}(x,t)\cos\rho t\,dt,
\]%
\[
S(\rho,x)=\frac{\sin\rho x}{\rho}+\int_{0}^{x}\mathbf{S}(x,t)\frac{\sin\rho
t}{\rho}dt
\]
for all $\rho\in\mathbb{C}$. These Volterra integral operators of the second
kind are known as the transformation operators. Among the properties of the
kernels $\mathbf{G}_{h}(x,t)$ and $\mathbf{S}(x,t)$ we will use the equalities%
\begin{equation}
\mathbf{G}_{h}(x,x)=h+\omega(x),\qquad\mathbf{S}(x,x)=\omega(x).
\label{G and S diag}%
\end{equation}

In \cite{KNT} the following series expansions for $\mathbf{G}_{h}(x,t)$ and
$\mathbf{S}(x,t)$ were obtained%
\[
\mathbf{G}_{h}(x,t)=%
{\displaystyle\sum\limits_{n=0}^{\infty}}
\frac{g_{n}(x)}{x}P_{2n}\left(  \frac{t}{x}\right)  ,\qquad\mathbf{S}\left(
x,t\right)  =\sum_{n=0}^{\infty}\frac{s_{n}(x)}{x}P_{2n+1}\left(  \frac{t}%
{x}\right)  ,
\]
where $P_{k}$ stands for the Legendre polynomial of degree $k$, and the
coefficients $g_{n}(x)$ and $s_{n}(x)$ are those from Theorem \ref{Th NSBF}.
Since $P_{k}(1)=1$ for all $k\in\mathbb{N}\cup\left\{  0\right\}  $, we have
\[%
{\displaystyle\sum\limits_{n=0}^{\infty}}
\frac{g_{n}(x)}{x}=h+\omega(x),\qquad\sum_{n=0}^{\infty}\frac{s_{n}(x)}%
{x}=\omega(x).
\]
Thus,
\begin{equation}
\frac{1}{x}%
{\displaystyle\sum\limits_{n=0}^{\infty}}
\left(  g_{n}(x)-s_{n}(x)\right)  =h. \label{sum = h}%
\end{equation}

\end{remark}

\begin{remark}
\label{Rem Recovery q, h, H}Since $\varphi_{h}(0,x)=g_{0}(x)+1$, we have that
\[
q(x)=\frac{\varphi_{h}^{\prime\prime}(0,x)}{\varphi_{h}(0,x)}=\frac
{g_{0}^{\prime\prime}(x)}{g_{0}(x)+1}.
\]
That is, $q(x)$ can be recovered directly from the first coefficient
$g_{0}(x)$, or, analogously, from $\psi_{0}(x)$:
\[
q(x)=\frac{\psi_{0}^{\prime\prime}(x)}{\psi_{0}(x)+1},
\]
or from a linear combination of $g_{0}(x)$ and $\psi_{0}(x)$. Moreover, we
have
\[
g_{0}^{\prime}(0)=\varphi_{h}^{\prime}(0,0)=h\quad\text{and}\quad\psi
_{0}^{\prime}(b)=\psi_{H}^{\prime}(0,b)=-H,
\]
and hence the constants $h$ and $H$ can also be recovered from the first NSBF
coefficients $g_{0}(x)$ and $\psi_{0}(x)$.
\end{remark}

\section{Approximate solution of inverse problem (IP1)\label{Sect4}}

The overall approach to solving the inverse problem (IP1) involves working
with the NSBF coefficients of the solutions and their derivatives. In the
first step, we solve a system of linear algebraic equations to compute the
characteristic functions $\Delta(\rho)$ and $\Delta^{0}(\rho)$, while in the
second, with their aid we construct another system of linear algebraic
equations for the NSBF coefficients of the solutions $\varphi_{h}(\rho,x)$,
$S(\rho,x)$ and $\psi_{H}(\rho,x)$. For this, identity (\ref{psiH=}) is used.
Finally the potential $q(x)$ and the constants $h$, $H$ are computed from the
first NSBF coefficients $g_{0}(x)$ and $\psi_{0}(x)$ as explained in Remark
\ref{Rem Recovery q, h, H}.

\subsection{Reconstruction of characteristic functions $\Delta(\rho)$ and
$\Delta^{0}(\rho)$}

\subsubsection{Recovering $\Delta^{0}(\rho)$}

For any singular number $\mu_{k}$ of problem $L^{0}$ we have the equality
$\psi_{H}(\mu_{k},0)=0$, which can be written in the form%
\[
\sum_{n=0}^{\infty}(-1)^{n}\psi_{n}(0)\mathbf{j}_{2n}(\mu_{k}b)=-\cos\left(
\mu_{k}b\right)  .
\]

Let us assume that a number of the eigenpairs $\rho_{k}=\sqrt{\lambda_{k}}$
and $\mu_{k}=\sqrt{\lambda_{k}^{0}}$ are given for $k=0,\ldots,K$. Considering
a partial sum of the series and writing down a corresponding equation for each
known $\mu_{k}$, we obtain the system of linear algebraic equations for the
coefficients $\psi_{n}(0)$:%
\begin{equation}
\sum_{n=0}^{N_{1}}(-1)^{n}\psi_{n}(0)\mathbf{j}_{2n}(\mu_{k}b)=-\cos\left(
\mu_{k}b\right)  ,\quad k=0,\ldots,K. \label{sys1}%
\end{equation}
Here $N_{1}$ needs to be less or equal $K$. Solving this system leads to an
approximation of the characteristic function $\Delta^{0}(\rho)$, which we
denote by $\Delta_{N_{1}}^{0}(\rho)$:
\[
\Delta_{N_{1}}^{0}(\rho)=\cos\left(  \rho b\right)  +\sum_{n=0}^{N_{1}%
}(-1)^{n}\psi_{n}(0)\mathbf{j}_{2n}(\rho b).
\]
It gives us an approximation of $\Delta^{0}(\rho)$. According to the estimates
from Theorem \ref{Th NSBF}, it is expected to provide more accurate results
for $\rho$ belonging to a strip $\left\vert \operatorname{Im}\rho\right\vert
\leq a$, where $a\geq0$ should not be too large.

\subsubsection{Recovering $\Delta(\rho)$ \label{subsubsect recovering Delta}}

Next, we consider the singular numbers $\rho_{k}$. We have the equality
$\Delta(\rho_{k})=0$, which can be written as $\varphi_{h}^{\prime}(\rho
_{k},b)+H\varphi_{h}(\rho_{k},b)=0$ or, equivalently,%
\[
\omega_{h,H}\cos\left(  \rho_{k}b\right)  +\sum_{n=0}^{\infty}(-1)^{n}%
\gamma_{n}(b)\mathbf{j}_{2n}(\rho_{k}b)+H\sum_{n=0}^{\infty}(-1)^{n}%
g_{n}(b)\mathbf{j}_{2n}(\rho_{k}b)=\rho_{k}\sin\left(  \rho_{k}b\right)  ,
\]
where $\omega_{h,H}:=h+H+\omega(b)$. Denote
\[
h_{n}:=(-1)^{n}\left(  \gamma_{n}(b)+Hg_{n}(b)\right)  ,\quad n=0,1,\ldots.
\]
Then
\[
\omega_{h,H}\cos\left(  \rho_{k}b\right)  +\sum_{n=0}^{\infty}h_{n}%
\mathbf{j}_{2n}(\rho_{k}b)=\rho_{k}\sin\left(  \rho_{k}b\right)  .
\]
Considering a partial sum of the series and writing down corresponding
equations for all available values $\rho_{k}$, we obtain a system of linear
algebraic equations for the unknowns $\omega_{h,H}$ and $\left\{
h_{n}\right\}  _{n=0}^{N_{2}}$ in the form%
\begin{equation}
\omega_{h,H}\cos\left(  \rho_{k}b\right)  +\sum_{n=0}^{N_{2}}h_{n}%
\mathbf{j}_{2n}(\rho_{k}b)=\rho_{k}\sin\left(  \rho_{k}b\right)  ,\quad
k=0,\ldots,K. \label{sys2}%
\end{equation}
Here $N_{2}\leq K-1$. For simplicity, we will assume that $N_{2}=N_{1}$, and
hence this is a system of linear algebraic equations for $N_{1}+2$ unknowns:
$\omega_{h,H}$, $\left\{  h_{n}\right\}  _{n=0}^{N_{1}}$. Solving (\ref{sys2})
gives us an approximation of $\Delta(\rho)$:%
\[
\Delta_{N_{2}}(\rho)=\omega_{h,H}\cos\left(  \rho b\right)  -\rho\sin\left(
\rho b\right)  +\sum_{n=0}^{N_{1}}h_{n}\mathbf{j}_{2n}(\rho b).
\]

\subsubsection{Choosing the number of the
coefficients\label{subsect choosing N1}}

For a detailed analysis of systems of the type (\ref{sys1}) and (\ref{sys2})
we refer to \cite{KT2021 IP1}, where it was explained that making (\ref{sys1})
and (\ref{sys2}) square sometimes can lead to worse numerical results, because
the systems may become ill-conditioned. However, no criterion was proposed for
choosing an optimal value of $N_{1}$. Here we fill this gap by proposing such
a simple and viable criterion. To this purpose, it is convenient to compute a
number of the coefficients $\left\{  s_{n}(b),\,g_{n}(b)\right\}
_{n=0}^{N_{1}}$. For this, note that since $\Delta^{0}(\mu_{k})=0$, from
(\ref{psiH=}) we have that $\psi_{H}(\mu_{k},x)=-\Delta\left(  \mu_{k}\right)
S(\mu_{k},x)$. Hence for $x=b$:
\begin{equation}
S(\mu_{k},b)=-\frac{1}{\Delta\left(  \mu_{k}\right)  }. \label{Smuk}%
\end{equation}
The values of $\Delta\left(  \mu_{k}\right)  $ are computed with the aid of
the computed previously $\omega_{h,H}$, $\left\{  h_{n}\right\}  _{n=0}%
^{N_{1}}$. We have $\Delta\left(  \mu_{k}\right)  \cong\Delta_{N_{1}}\left(
\mu_{k}\right)  =\omega_{h,H}\cos\left(  \mu_{k}b\right)  -\mu_{k}\sin\left(
\mu_{k}b\right)  +\sum_{n=0}^{N_{1}}h_{n}\mathbf{j}_{2n}(\mu_{k}b)$. Thus, for
the coefficients $\left\{  s_{n}(b)\right\}  _{n=0}^{N_{1}}$ we obtain the
system%
\begin{equation}
\sum_{n=0}^{N_{1}}(-1)^{n}s_{n}(b)\mathbf{j}_{2n+1}(\mu_{k}b)=-\left(  \sin
\mu_{k}b+\frac{\mu_{k}}{\Delta_{N_{1}}\left(  \mu_{k}\right)  }\right)  ,\quad
k=0,\ldots,K. \label{sys3}%
\end{equation}
Analogously, observe that for $\rho=\rho_{k}$ identity (\ref{psiH=}) gives
$\psi_{H}(\rho_{k},x)=\Delta^{0}\left(  \rho_{k}\right)  \varphi_{h}(\rho
_{k},x)$, which for $x=b$ leads to the equality%
\[
\varphi_{h}(\rho_{k},x)=\frac{1}{\Delta^{0}\left(  \rho_{k}\right)  }.
\]
Thus, a number of the coefficients $\left\{  g_{n}(b)\right\}  _{n=0}^{N_{1}}$
can be computed from the system
\begin{equation}
\sum_{n=0}^{N_{1}}(-1)^{n}g_{n}(b)\mathbf{j}_{2n}(\rho_{k}b)=\frac{1}%
{\Delta_{N_{1}}^{0}\left(  \rho_{k}\right)  }-\cos\rho_{k}b,\quad
k=0,\ldots,K, \label{sys4}%
\end{equation}
where the values of $\Delta_{N_{1}}^{0}\left(  \rho_{k}\right)  $ are computed
with the aid of the computed previously $\left\{  \psi_{n}(0)\right\}
_{n=0}^{N_{1}}$: $\Delta_{N_{1}}^{0}\left(  \rho_{k}\right)  =\cos\left(
\rho_{k}b\right)  +\sum_{n=0}^{N_{1}}(-1)^{n}\psi_{n}(0)\mathbf{j}_{2n}%
(\rho_{k}b)$.

Now, for a given $K$, a convenient choice of $N_{1}$ can be made in two ways,
both by using identity (\ref{psiH=}). For $x=b$, identity (\ref{psiH=}) turns
into the equality%
\[
\Delta^{0}\left(  \rho\right)  \varphi_{h}(\rho,b)-\Delta\left(  \rho\right)
S(\rho,b)=1.
\]
Choosing a set of points $\rho=r_{j}$, $j=1,\ldots,J$, where $J$ can be
arbitrarily large, and $\left\vert \operatorname{Im}r_{j}\right\vert \leq a$
(with a reasonably large $a$ or $a=0$), one can compute the functional
\[
P(N_{1}):=\max_{j}\left\vert \Delta_{N_{1}}^{0}\left(  r_{j}\right)
\varphi_{h,N_{1}}(r_{j},b)-\Delta_{N_{1}}\left(  r_{j}\right)  S_{N_{1}}%
(r_{j},b)-1\right\vert ,
\]
which provides information on the accuracy attained for the chosen $N_{1}$.
Minimization of $P(N_{1})$ gives us a suitable choice of $N_{1}$.

Another criterion for choosing $N_{1}$ can be devised as follows. Consider
identity (\ref{psiH=}) for $x=b$ and $\rho=0$:%
\[
1=\left(  1+\psi_{0}(0)\right)  \left(  1+g_{0}(b)\right)  -b\left(
\omega_{h,H}+h_{0}\right)  \left(  1+\frac{s_{0}(b)}{3}\right)  ,
\]
or equivalently,
\[
g_{0}(b)\left(  1+\psi_{0}(0)\right)  +\psi_{0}(0)-\frac{b}{3}\left(
\omega_{h,H}+h_{0}\right)  \left(  3+s_{0}(b)\right)  =0.
\]
Define
\[
R(N_{1}):=\left\vert g_{0}(b)\left(  1+\psi_{0}(0)\right)  +\psi_{0}%
(0)-\frac{b}{3}\left(  \omega_{h,H}+h_{0}\right)  \left(  3+s_{0}(b)\right)
\right\vert ,
\]
where all the magnitudes are computed with the parameter $N_{1}$ in systems
(\ref{sys1}), (\ref{sys2}), (\ref{sys3}), (\ref{sys4}). Minimizing $R(N_{1})$
gives us another criterion for the choice of a suitable value of $N_{1}$.
Obviously, this criterion is especially useful when zero is not an eigenvalue
neither of problem $L$ nor of $L^{0}$.

It turns out that both criteria deliver similar results, as we show below.
Thus, one can use the value $N_{1}$ which minimizes $R(N_{1})$, since the
computation of $R(N_{1})$ is easier and faster.

\subsection{Reconstruction of $q(x)$, $h$ and $H$
\label{subsect reconstruction q and h}}

In the second step we use identity (\ref{psiH=}), where $\Delta^{0}\left(
\rho\right)  $ and $\Delta\left(  \rho\right)  $ are already approximated by
$\Delta_{N_{1}}^{0}\left(  \rho\right)  $ and $\Delta_{N_{1}}\left(
\rho\right)  $, respectively, and the solutions $\psi_{H}(\rho,x)$,
$\varphi_{h}(\rho,x)$ and $S(\rho,x)$ are replaced by their NSBF
representations. Thus, choosing a set of points $\left\{  r_{j}\right\}
_{j=1}^{J}\subset\mathbb{R}$ sufficiently large, for any $x\in\left(
0,b\right)  $ we consider the system of linear algebraic equations for the
NSBF coefficients%
\begin{gather}
\Delta_{N_{1}}^{0}\left(  r_{j}\right)  \sum_{n=0}^{N_{3}}(-1)^{n}%
g_{n}(x)\mathbf{j}_{2n}(r_{j}x)-\frac{\Delta\left(  r_{j}\right)  }{r_{j}}%
\sum_{n=0}^{N_{3}}(-1)^{n}s_{n}(x)\mathbf{j}_{2n+1}(r_{j}x)-\sum_{n=0}^{N_{3}%
}(-1)^{n}\psi_{n}(x)\mathbf{j}_{2n}(r_{j}\left(  x-b\right)  )\nonumber\\
=\cos\left(  r_{j}\left(  x-b\right)  \right)  -\Delta_{N_{1}}^{0}\left(
r_{j}\right)  \cos\left(  r_{j}x\right)  +\frac{\Delta\left(  r_{j}\right)
}{r_{j}}\sin\left(  r_{j}x\right)  ,\quad j=1,\ldots,J, \label{sys main1}%
\end{gather}
where $3(N_{3}+1)\leq J$. Solving this system for a sufficiently dense set of
points $x\in\left(  0,b\right)  $ gives us the values of $g_{0}(x)$ and
$\psi_{0}(x)$. Note that there is no need to solve system (\ref{sys main1}) at
the endpoints, because the values of $g_{0}(x)$ and $\psi_{0}(x)$ there are
already known. We have that $g_{0}(0)=\psi_{0}(b)=0$ (see (\ref{g0s0}) and
(\ref{psi0t0})), while the values\ $g_{0}(b)$ and $\psi_{0}(0)$ were computed
in the first step.

Finally, $q(x)$, $h$ and $H$ are recovered from $g_{0}(x)$ and $\psi_{0}(x)$
as explained in Remark \ref{Rem Recovery q, h, H}.

\section{Approximate solution of other inverse spectral
problems\label{Sect Other inverse problems}}

In this section we show that the approach presented above is quite universal
and can be easily adapted for solving a wide variety of inverse spectral
problems. In particular, we consider the recovery from Weyl's function; from
one spectrum and multiplier constants (which are values of normalized
eigenfunctions at endpoints); from one spectrum and a sequence of norming constants.

\subsection{Recovery from Weyl function}

Let $\Phi(\rho,x)$ be a solution of (\ref{Schr}) satisfying the boundary
conditions%
\[
U(\Phi)=1,\quad V(\Phi)=0.
\]
The solution $\Phi(\rho,x)$ is known as the Weyl solution. The Weyl function
of problem $L$ is then defined as
\[
M(\rho):=\Phi(\rho,0)
\]
(see, e.g., \cite[Sect. 1.2]{Yurko2007}).

Consider the following inverse problem.

\textbf{Inverse Problem (IP2) }Given the Weyl function $M(\rho)$, find
$q(x)\in%
\mathcal{L}%
^{2}\left(  0,b\right)  $ and the constants $h$, $H$, such that $M(\rho)$ be
the Weyl function of problem $L$.

More precisely, we are interested in the approximate solution of (IP2) when
$M(\rho)$ is given on a finite set of points $\left\{  \rho=z_{k}\right\}
_{k=1}^{K_{1}}$.

It is easy to see that
\[
\Phi(\rho,x)=S(\rho,x)+M(\rho)\varphi_{h}(\rho,x).
\]
Hence
\[
V(\Phi)=S^{\prime}(\rho,b)+HS(\rho,b)+M(\rho)\left(  \varphi_{h}^{\prime}%
(\rho,b)+H\varphi_{h}(\rho,b)\right)  =0,
\]
or equivalently,%
\begin{equation}
\Delta^{0}\left(  \rho\right)  +M(\rho)\Delta\left(  \rho\right)  =0.
\label{identity Deltas M}%
\end{equation}
Substitution of $\rho=z_{k}$, $k=1,\ldots,K_{1}$ and of the NSBF
representations (truncated up to some $N_{1}$) leads to the system of linear
algebraic equations for the NSBF coefficients $\omega_{h,H}$, $\left\{
\psi_{n}(0)\right\}  $ and $\left\{  h_{n}\right\}  $:
\begin{equation}
\sum_{n=0}^{N_{1}}(-1)^{n}\psi_{n}(0)\mathbf{j}_{2n}(z_{k}b)+M(z_{k}%
)\omega_{h,H}\cos\left(  z_{k}b\right)  +M(z_{k})\sum_{n=0}^{N_{1}}%
h_{n}\mathbf{j}_{2n}(z_{k}b)=M(z_{k})z_{k}\sin\left(  z_{k}b\right)
-\cos\left(  z_{k}b\right)  , \label{sys Weyl}%
\end{equation}
$k=1,\ldots,K_{1}$, $2(N_{1}+1)\leq K_{1}-1$. Solving this system gives us the
approximations $\Delta_{N_{1}}^{0}\left(  \rho\right)  $ and $\Delta_{N_{1}%
}\left(  \rho\right)  $, after which we proceed to the second step, see
subsection \ref{subsect reconstruction q and h}.

\subsection{Recovery from eigenvalues and multiplier constants}

Let $\rho_{k}$ be a singular value of problem $L$. We have that the solutions
$\varphi_{h}(\rho_{k},x)$ and $\psi_{H}(\rho_{k},x)$ are necessarily linearly
dependent, so there exists a constant $\beta_{k}\neq0$ such that
\begin{equation}
\varphi_{h}(\rho_{k},x)=\beta_{k}\psi_{H}(\rho_{k},x). \label{phih=psiH}%
\end{equation}
Choosing here $x=0$ and $x=b$, we obtain that
\begin{equation}
\beta_{k}=\frac{1}{\psi_{H}(\rho_{k},0)}=\varphi_{h}(\rho_{k},b).
\label{beta_k=}%
\end{equation}
Now the following inverse problem is considered.

\textbf{Inverse Problem (IP3) }Given the singular values $\left\{  \rho
_{k}\right\}  _{k=0}^{\infty}$ and the multiplier constants $\left\{
\beta_{k}\right\}  _{k=0}^{\infty}$ of problem $L$, find $q(x)\in%
\mathcal{L}%
^{2}\left(  0,b\right)  $ and the constants $h$, $H$.

A similar inverse problem was considered in \cite{Brown et al 2003}, however
with boundary conditions not containing unknown constants and for real valued
$q(x)$. We also note that it is often convenient to reduce other types of
inverse spectral problems to (IP3) (see \cite{AKK2024IPI}, \cite{AKK2024MMAS},
\cite{CKK2024MMAS}, \cite{KrSpectrumCompl}). Our approach allows us to solve
this problem numerically in one step. Indeed, equality (\ref{phih=psiH}) leads
to the following system of linear algebraic equations for the NSBF
coefficients of the solutions $\varphi_{h}(\rho,x)$ and $\psi_{H}(\rho,x)$:%
\begin{equation}
\sum_{n=0}^{\infty}(-1)^{n}g_{n}(x)\mathbf{j}_{2n}(\rho_{k}x)-\beta_{k}%
\sum_{n=0}^{\infty}(-1)^{n}\psi_{n}(x)\mathbf{j}_{2n}(\rho_{k}\left(
x-b\right)  )=\beta_{k}\cos\left(  \rho_{k}\left(  x-b\right)  \right)
-\cos\left(  \rho_{k}x\right)  , \label{sysIP3}%
\end{equation}
$k=0,1,\ldots$. Solving this system, or its truncated version approximately on
a sufficiently dense set of points $x\in\left[  0,b\right]  $, we find
$g_{0}(x)$ and $\psi_{0}(x)$ and apply Remark \ref{Rem Recovery q, h, H} for
recovering $q(x)$, $h$ and $H$. Here unlike the solution of (\ref{sys main1})
in subsection \ref{subsect reconstruction q and h}, we solve (\ref{sysIP3}) at
the endpoints as well, where it is convenient to take into account that
$g_{n}(0)=\psi_{n}(b)=0$, $n=0,1,\ldots$.

This direct way to solve (IP3) in one step, possibly, is not the best option
for numerical computations. It requires a relatively large number of the data
to be given. Another possibility leads to overall better numerical results.
First, by the given $\rho_{k}$ we approximate $\Delta\left(  \rho\right)  $.
This is done exactly as in subsection \ref{subsubsect recovering Delta}. Next,
we use $\beta_{k}$ for computing a number of the coefficients $\psi_{n}(0)$.
Namely, we have $\psi_{H}(\rho_{k},0)=1/\beta_{k}$, $k=0,1,\ldots$. This gives
us the system of linear algebraic equations for the coefficients $\psi_{n}%
(0)$:%
\[
\sum_{n=0}^{N_{1}}(-1)^{n}\psi_{n}(0)\mathbf{j}_{2n}(\rho_{k}b)=\frac{1}%
{\beta_{k}}-\cos\left(  \rho_{k}b\right)  ,\quad k=0,\ldots,K
\]
(compare to (\ref{sys1})).

Thus, again, we have approximated $\Delta^{0}\left(  \rho\right)  $ and
$\Delta\left(  \rho\right)  $ by $\Delta_{N_{1}}^{0}\left(  \rho\right)  $ and
$\Delta_{N_{1}}\left(  \rho\right)  $, which leads us to the second step of
the method, explained in subsection \ref{subsect reconstruction q and h}.

\subsection{Recovery from eigenvalues and norming constants
\label{Subsect IP4}}

For every singular value $\rho_{k}$ of problem $L$, we denote by $\alpha_{k}$
the corresponding norming constant%
\[
\alpha_{k}:=\int_{0}^{b}\varphi_{h}^{2}(\rho_{k},x)dx.
\]

Consider the inverse problem.

\textbf{Inverse Problem (IP4) }Given the singular values $\left\{  \rho
_{k}\right\}  _{k=0}^{\infty}$ and the norming constants $\left\{
\alpha\right\}  _{k=0}^{\infty}$ of problem $L$, find $q(x)\in%
\mathcal{L}%
^{2}\left(  0,b\right)  $ and the constants $h$, $H$.

We recall that
\begin{equation}
\frac{\alpha_{k}}{\beta_{k}}=-\frac{\overset{\cdot}{\Delta}(\rho_{k})}%
{2\rho_{k}}, \label{equality alpha beta k}%
\end{equation}
where $\overset{\cdot}{\Delta}$ means the derivative with respect to $\rho$.
Indeed, with the aid of the equality
\[
\frac{d}{dx}W\left[  \psi_{H}(\rho,x),\varphi_{h}(\rho_{k},x)\right]  =\left(
\rho^{2}-\rho_{k}^{2}\right)  \psi_{H}(\rho,x)\varphi_{h}(\rho_{k},x)
\]
we obtain%
\[
\left(  \rho^{2}-\rho_{k}^{2}\right)  \int_{0}^{b}\psi_{H}(\rho,x)\varphi
_{h}(\rho_{k},x)dx=\left.  W\left[  \psi_{H}(\rho,x),\varphi_{h}(\rho
_{k},x)\right]  \right\vert _{0}^{b}=-\Delta\left(  \rho\right)  ,
\]
or
\[
\frac{\Delta\left(  \rho\right)  }{\left(  \lambda-\lambda_{k}\right)  }%
=-\int_{0}^{b}\psi_{H}(\rho,x)\varphi_{h}(\rho_{k},x)dx,
\]
which at the limit $\lambda\rightarrow\lambda_{k}$ gives us
(\ref{equality alpha beta k}).

Thus,
\begin{equation}
\beta_{k}=-\frac{2\rho_{k}\alpha_{k}}{\overset{\cdot}{\Delta}(\rho_{k})}.
\label{beta k =}%
\end{equation}
It is easy to obtain an NSBF representation for $\overset{\cdot}{\Delta}%
(\rho)$. We have
\[
\Delta\left(  \rho\right)  =\omega_{h,H}\cos\left(  \rho b\right)  -\rho
\sin\left(  \rho b\right)  +\sum_{n=0}^{\infty}h_{n}\mathbf{j}_{2n}(\rho b).
\]
Hence
\begin{equation}
\overset{\cdot}{\Delta}(\rho)=-\left(  1+b\omega_{h,H}\right)  \sin\left(
\rho b\right)  -\rho b\cos\left(  \rho b\right)  +\sum_{n=0}^{\infty}%
h_{n}\left(  b\mathbf{j}_{2n+1}(\rho b)-\frac{2n}{\rho}\mathbf{j}_{2n}(\rho
b)\right)  . \label{Delta dot}%
\end{equation}
Now, the algorithm for solving (IP4) can be summarized as follows. First, from
(\ref{sys2}) we compute $\omega_{h,H}$, $\left\{  h_{n}\right\}  _{n=0}%
^{N_{2}}$. Next, with the aid of these coefficients the values $\overset
{\cdot}{\Delta}(\rho_{k})$ are computed approximately from (\ref{Delta dot}),%
\[
\overset{\cdot}{\Delta}(\rho_{k})\cong\overset{\cdot}{\Delta}_{N_{2}}(\rho
_{k})=-\left(  1+b\omega_{h,H}\right)  \sin\left(  \rho_{k}b\right)  -\rho
_{k}b\cos\left(  \rho_{k}b\right)  +\sum_{n=0}^{N_{2}}h_{n}\left(
b\mathbf{j}_{2n+1}(\rho_{k}b)-\frac{2n}{\rho_{k}}\mathbf{j}_{2n}(\rho
_{k}b)\right)  .
\]
Now, the multiplier constants $\beta_{k}$ are computed from (\ref{beta k =}).
Thus, problem (IP4) is reduced to problem (IP3), and one can apply the
algorithm described in the previous section.

\section{Numerical results\label{Sect Numeric}}

The proposed approach can be implemented directly using an available numerical
computing environment. All the reported computations were performed in Matlab
R2024a on an Intel i7-1360P equipped laptop computer and took no more than
several seconds.

\subsection{Solution of two-spectrum inverse problem (IP1)}

\textbf{Example 1. }Consider equation (\ref{Schr}) with
\begin{equation}
q(x)=x^{2},\quad0<x<1, \label{q Ex1}%
\end{equation}
and the boundary conditions (\ref{bc1}) with $h=10$ and $H=\pi$. For the first
test we recover this Sturm-Liouville problem from ten eigenpairs: $\left\{
\lambda_{k},\,\lambda_{k}^{0}\right\}  _{k=0}^{9}$. The \textquotedblleft
exact\textquotedblright\ eigenvalues were computed by Matslise \cite{Ledoux
Matslise}. The minima of both functionals $P(N_{1})$ and $R(N_{1})$ (see
subsection \ref{subsect choosing N1}) were attained at $N_{1}=7$, see Fig. 1.
Thus, in the first step both systems (\ref{sys1}) and (\ref{sys2}) were
considered with $K=9$ and $N_{1}=N_{2}=7$. It is worth noticing that the
parameter $\omega_{h,H}\cong13.3082593201899$ was computed from system
(\ref{sys2}) with an absolute error $3.6\cdot10^{-6}$.
\begin{figure}
[ptb]
\begin{center}
\includegraphics[
height=2.7826in,
width=3.7077in
]%
{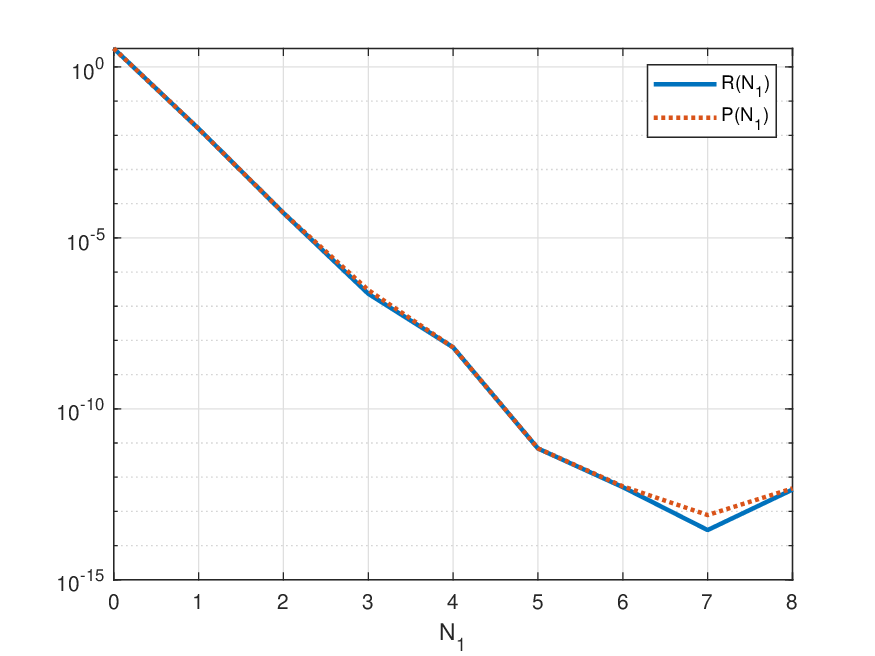}%
\caption{$P(N_{1})$ and $R(N_{1})$ computed for Example 1, with $10$
eigenpairs given. Both indicate the choice of $N_{1}=7$.}%
\label{Fig1}%
\end{center}
\end{figure}
For the second step, with the aid of the computed values $\omega_{h,H}$,
$\left\{  \psi_{n}(0),\,h_{n}\right\}  _{n=0}^{7}$, the approximate
characteristic functions $\Delta_{7}^{0}\left(  \rho\right)  $ and $\Delta
_{7}\left(  \rho\right)  $ were computed at $J=1501$ points $\rho=r_{j}$
chosen as follows $r_{j}=10^{\alpha_{j}}$ with $\alpha_{j}$, $j=1,\ldots,J$,
being distributed uniformly on $\left[  \log(0.01),\log(1000)\right]  =\left[
-2,\,3\right]  $. This is a logarithmically spaced point distribution on the
segment $[0.01,1000]$, with points more densely distributed near $\rho=0.01$
and less densely distributed as $\rho$ increases. The parameter $N_{3}$ here
and in the other numerical tests was chosen equal to $N_{1}$. Thus, system
(\ref{sys main1}) consisted of $1501$ equations for $24$ unknowns. The result
of the reconstruction of $q(x)$ is shown in Fig. 2.%

\begin{figure}
[ptb]
\begin{center}
\includegraphics[
height=4.6009in,
width=6.1308in
]%
{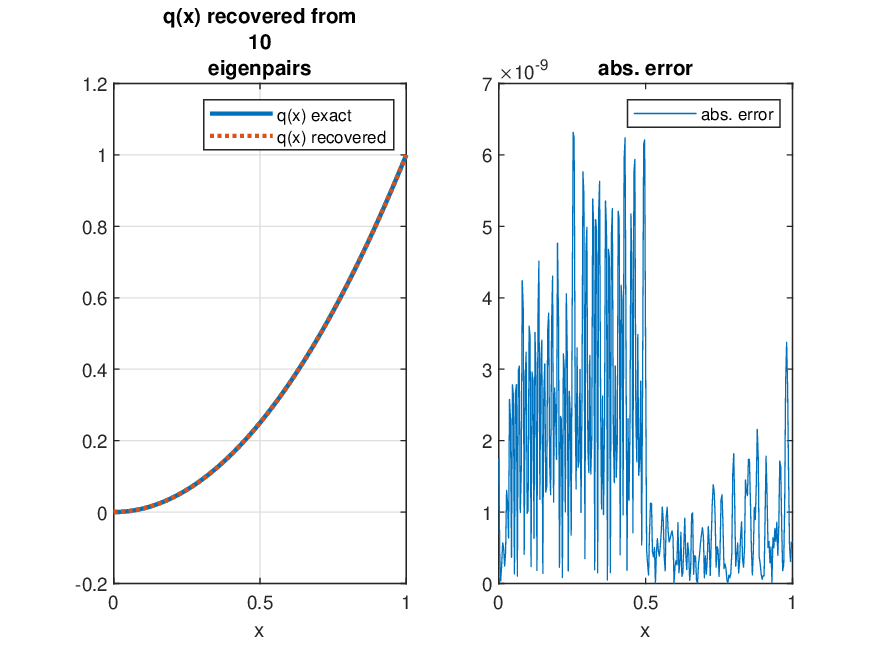}%
\caption{Potential from Example 1, recovered from $10$ eigenpairs. The maximum
absolute error of the recovered potential is $6.3\cdot10^{-9}$. The constants
$h=10$ and $H=\pi$ are recovered with the absolute error $3.4\cdot10^{-13}$
and $2.9\cdot10^{-12}$, respectively.}%
\label{Fig2}%
\end{center}
\end{figure}

The potential (\ref{q Ex1}) on the interval $\left(  0,1\right)  $ is a
relatively simple example and can be recovered from less eigenpairs. For
example, the maximum absolute error of its recovery from five eigenpairs was
$6.8\cdot10^{-3}$. The constants $h$ and $H$ were recovered with the absolute
error $4.7\cdot10^{-6}$ and $8.1\cdot10^{-6}$, respectively.

Moreover, the reconstruction procedure results to be stable with respect to a
certain level of noise in the input data. The noise is introduced as follows%
\begin{equation}
\lambda_{k,\text{noisy}}=\lambda_{k}+\sigma\sin(\frac{\left(  k+1\right)  \pi
}{37}),\quad k=0,\ldots,K \label{lambda noisy}%
\end{equation}
and similarly for $\lambda_{k}^{0}$. Here the coefficient $\sigma$ determines
the level of noise. In the considered example of the reconstruction from five
eigenpairs ($K=4$), we introduced the noise with $\sigma=0.001$. As a result,
the potential (\ref{q Ex1}) was still recovered quite accurately, with the
maximum absolute error of $0.047$. The constants $h$ and $H$ were recovered
with the absolute error $9.9\cdot10^{-5}$ and $3.1\cdot10^{-4}$, respectively.

In the rest of the numerical tests, for the system of equations
(\ref{sys main1}) arising in the last step of the method, we choose the same
distribution of the points $\rho=r_{j}$ , $j=1,\ldots,J$, as described above
in Example 1.

\textbf{Example 2. }Consider the potential
\[
q(x)=e^{x},\quad0<x<\pi,
\]
and the boundary conditions (\ref{bc1}) with $h=10$ and $H=\pi$. Sometimes
this potential is called the first Paine potential \cite{Paine}. The result of
its recovery from 15 eigenpairs is shown in Fig. 3. The maximum absolute error
resulted in $1.7\cdot10^{-4}$. The constants $h$ and $H$ were recovered with
the absolute error\ $1.6\cdot10^{-8}$ and $6.2\cdot10^{-7}$, respectively. The
parameter $N_{1}$ ($=N_{2}$), as a minimum of $R(N_{1}),$ resulted in its
maximum possible value $N_{1}=N_{2}=13$, which corresponds to the case when
system (\ref{sys2}) is square.%

\begin{figure}
[ptb]
\begin{center}
\includegraphics[
height=5.0433in,
width=6.7211in
]%
{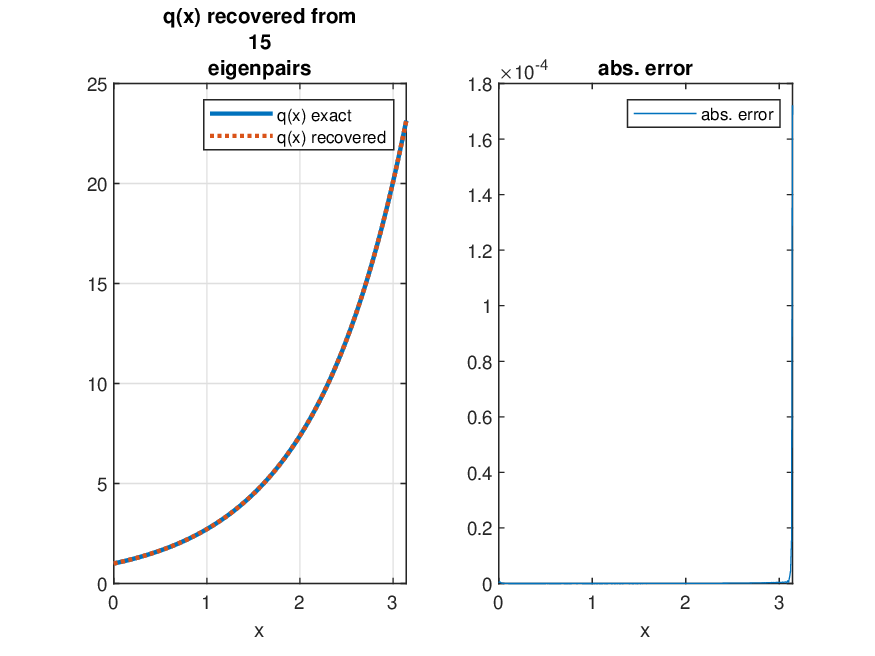}%
\caption{First Paine potential (Example 2) recovered from $15$ eigenpairs. The
maximum absolute error is $1.7\cdot10^{-4}$. The constants $h$ and $H$ are
recovered with the absolute error\ $1.6\cdot10^{-8}$ and $6.2\cdot10^{-7}$,
respectively.}%
\label{Fig3}%
\end{center}
\end{figure}

Next, we recovered the same Sturm-Liouville problem from $15$ noisy eigenpairs
with $\sigma=0.01$. This level of noise is high. However, the problem was
recovered sufficiently accurately. Fig. 4 shows the result of the recovery of
the potential. The maximum absolute error attained at the endpoint $x=\pi$
resulted in $0.34$. The constants $h$ and $H$ were recovered with the absolute
error\ $2.07\cdot10^{-4}$ and $0.02$, respectively.%

\begin{figure}
[ptb]
\begin{center}
\includegraphics[
height=4.8308in,
width=6.4383in
]%
{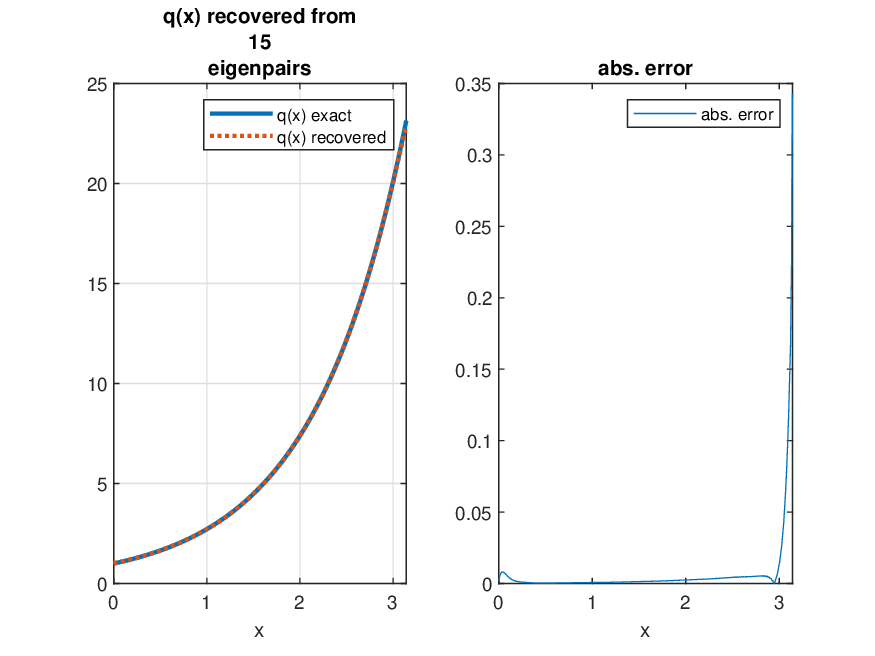}%
\caption{Same as Fig. 3, but recovered from $15$ noisy eigenpairs with a high
level of noise $\sigma=0.01$ (see (\ref{lambda noisy})). The maximum absolute
error attained at the endpoint $x=\pi$ resulted in $0.34$. The constants $h$
and $H$ were recovered with the absolute error\ $2.07\cdot10^{-4}$ and $0.02$,
respectively.}%
\label{Fig4}%
\end{center}
\end{figure}

The method works equally well for complex-valued potentials. Consider now
\begin{equation}
q(x)=e^{x}+\pi i,\quad0<x<\pi, \label{Paine1+Im}%
\end{equation}
and the same constants $h$ and $H$. Obviously, since the imaginary part of
$q(x)$ is constant, the eigenvalues of the corresponding Sturm-Liouville
problems are obtained by adding $\pi i$ to the eigenvalues of the real valued
potential $\operatorname{Re}q(x)$. The result of the recovery of $q(x)$ from
$15$ noisy eigenpairs with $\sigma=0.01$, is presented in Fig. 5. The maximum
absolute error attained at the endpoint $x=\pi$ resulted in $1.1$. The
constants $h$ and $H$ were recovered with the absolute error\ $0.0056$ and
$0.0023$, respectively.%

\begin{figure}
[ptb]
\begin{center}
\includegraphics[
height=5.0853in,
width=6.7786in
]%
{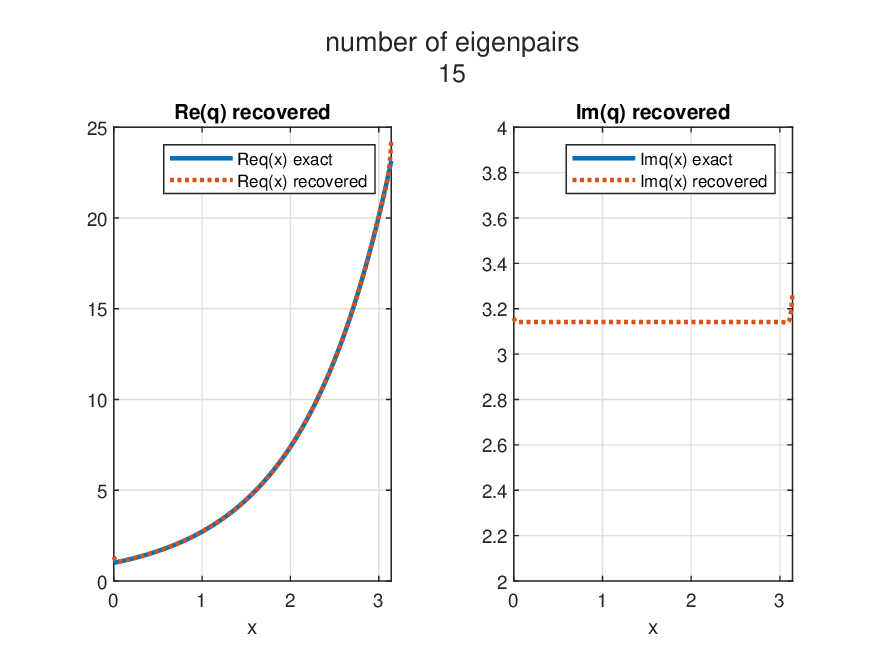}%
\caption{Complex valued otential (\ref{Paine1+Im}) recovered from $15$ noisy
eigenpairs with $\sigma=0.01$. The maximum absolute error attained at the
endpoint $x=\pi$ resulted in $1.1$. The constants $h$ and $H$ were recovered
with the absolute error\ $0.0056$ and $0.0023$, respectively.}%
\label{Fig5}%
\end{center}
\end{figure}

\textbf{Example 3. }Consider the complex-valued potential
\[
q(x)=\left(  x^{\frac{\pi}{2}}+\pi\right)  \cos\left(  8x\right)  +\pi
^{2}-i\sqrt{5},\quad0<x<\pi,
\]
and the boundary conditions (\ref{bc1}) with $h=\sqrt{2}$ and $H=-e$. Since
the imaginary part of $q(x)$ is constant, it is not difficult to obtain the
eigenpairs by computing them with Matslise for the real part and then adding
to them $\operatorname{Im}q(x)i$. The recovered potential is shown in Fig. 6.
It was computed from $50$ eigenpairs. The maximum absolute error of the
recovered potential resulted in $0.05$, attained at $x=0$. The parameter
$N_{1}$ ($=N_{2}$), as a minimum of $R(N_{1})$, resulted in $N_{1}=31$. The
constants $h$ and $H$ were recovered with the absolute error\ $2.8\cdot
10^{-4}$ and $1.1\cdot10^{-5}$, respectively.%

\begin{figure}
[ptb]
\begin{center}
\includegraphics[
height=4.7742in,
width=6.3625in
]%
{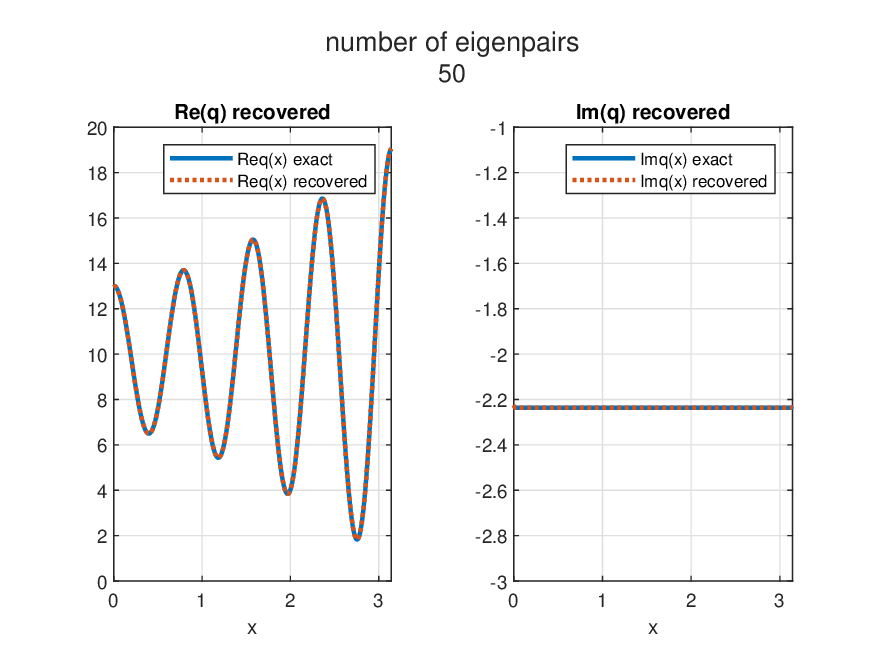}%
\caption{Complex valued potential from Example 3, recovered from $50$
eigenpairs. The maximum absolute error of the recovered potential resulted in
$0.05$, attained at $x=0$. The parameter $N_{1}$ ($=N_{2}$), as a minimum of
$R(N_{1})$, resulted in $N_{1}=31$. The constants $h=\sqrt{2}$ and $H=-e$ were
recovered with the absolute error\ $2.8\cdot10^{-4}$ and $1.1\cdot10^{-5}$,
respectively.}%
\label{Fig6}%
\end{center}
\end{figure}

Next, we reconstructed the same Sturm-Liouville problem from $50$ noisy
eigenpairs with $\sigma=0.001$. Fig. 7 shows the result. The maximum absolute
error of the recovered potential resulted in $0.7$, attained at $x=\pi$.
Again, as in the previous example, the noise in the input data resulted in a
smaller value of the parameter $N_{1}$, computed as a minimum of $R(N_{1})$.
We obtained $N_{1}=N_{2}=22$. The constants $h$ and $H$ were recovered with
the absolute error\ $6\cdot10^{-4}$ and $1.8\cdot10^{-3}$, respectively.%

\begin{figure}
[ptb]
\begin{center}
\includegraphics[
height=5.1054in,
width=6.8041in
]%
{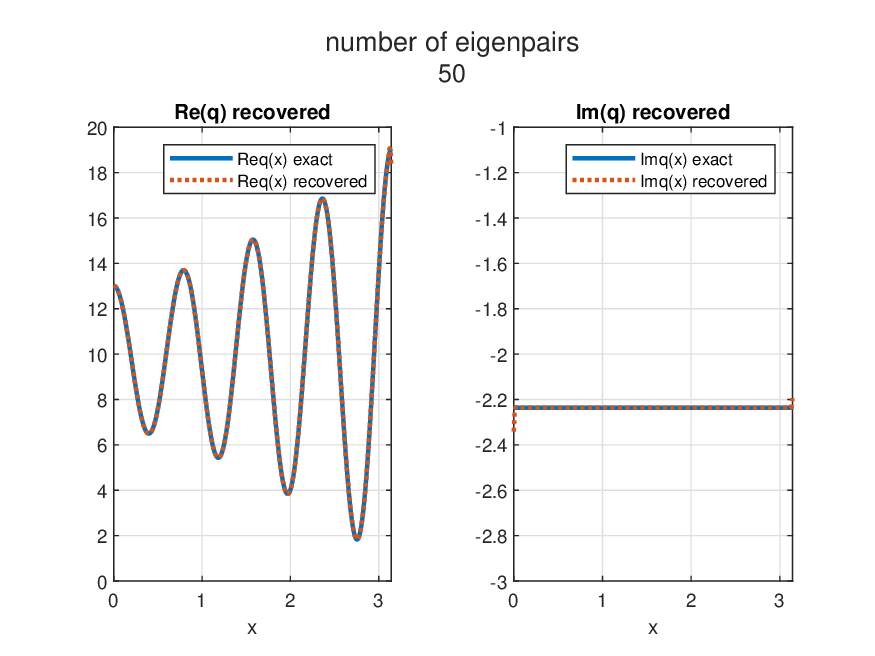}%
\caption{Complex valued potential from Example 3, recovered from $50$ noisy
eigenpairs with $\sigma=0.001$. The maximum absolute error of the recovered
potential resulted in $0.7$, attained at $x=\pi$. Again, as in the previous
example, the noise in the input data resulted in a smaller value of the
parameter $N_{1}$, computed as a minimum of $R(N_{1})$. We obtained
$N_{1}=N_{2}=22$. The constants $h$ and $H$ were recovered with the absolute
error\ $6\cdot10^{-4}$ and $1.8\cdot10^{-3}$, respectively.}%
\label{Fig7}%
\end{center}
\end{figure}

The computation of the eigenvalues for non-selfadjoint Sturm-Liouville
problems is, in general, a distinct and challenging task. In the previous two
numerical tests, the eigenvalues were computed for selfadjoint problems and
subsequently converted to those of the non-selfadjoint problems, as the
imaginary part of the potential was simply a constant. The next example is
different. It is purely imaginary, and the eigenvalues were computed in the
framework of \cite{KM Online MMAS} with the aid of the NSBF representations
and argument principle.

\textbf{Example 4. }Consider the Mathieu potential\textbf{ }
\begin{equation}
q(x)=is\cos(2x),\quad0<x<\pi, \label{Mathieu}%
\end{equation}
and the boundary conditions (\ref{bc1}) with $h=0.7$ and $H=i$. Let $s=2$.
Thus, this inverse problem involves a complex-valued potential and a complex
constant in the boundary condition. Fig. 8 shows the result of the recovery of
the potential (\ref{Mathieu}) from ten eigenpairs. The maximum absolute error
resulted in $3.4\cdot10^{-3}$. The constants $h$ and $H$ were recovered with
the absolute error $8.5\cdot10^{-5}$ and $1.2\cdot10^{-4}$, respectively.
Thus, these results are more accurate than those reported in \cite{KM Online
MMAS}.%

\begin{figure}
[ptb]
\begin{center}
\includegraphics[
height=5.2924in,
width=7.0541in
]%
{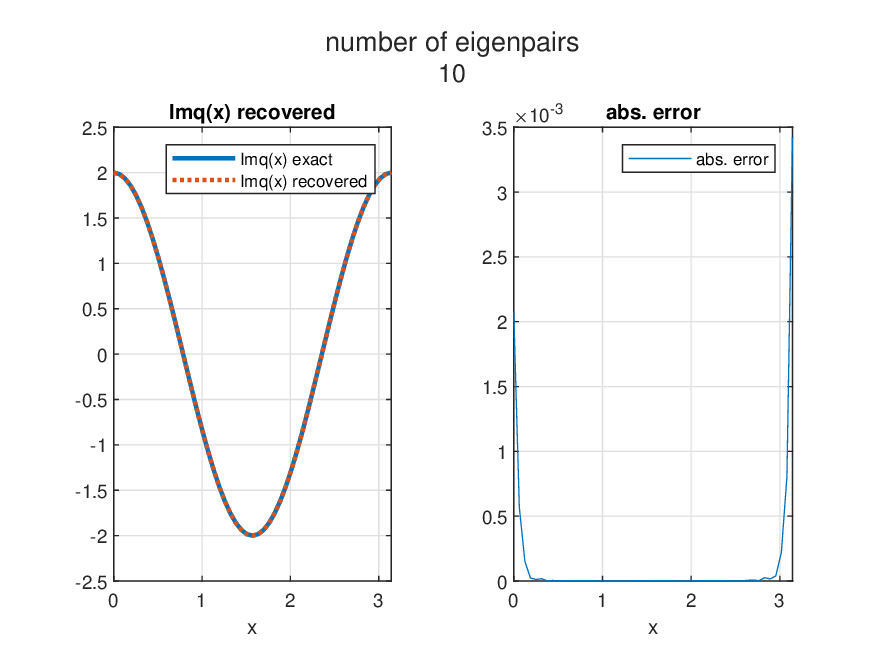}%
\caption{Mathieu potential from Example 4, recovered from ten eigenpairs. The
maximum absolute error resulted in $3.4\cdot10^{-3}$. The constants $h=0.7$
and $H=i$ were recovered with the absolute error $8.5\cdot10^{-5}$ and
$1.2\cdot10^{-4}$, respectively.}%
\label{Fig8}%
\end{center}
\end{figure}
Moreover, the method shows a remarkable stability. Fig. 9 presents the result
of the reconstruction of the potential from ten noisy eigenpairs with
$\sigma=0.001$. The maximum absolute error resulted in $0.043$. The constants
$h$ and $H$ were recovered with the absolute error $6\cdot10^{-5}$ and
$3\cdot10^{-3}$, respectively.%

\begin{figure}
[ptb]
\begin{center}
\includegraphics[
height=4.9028in,
width=6.535in
]%
{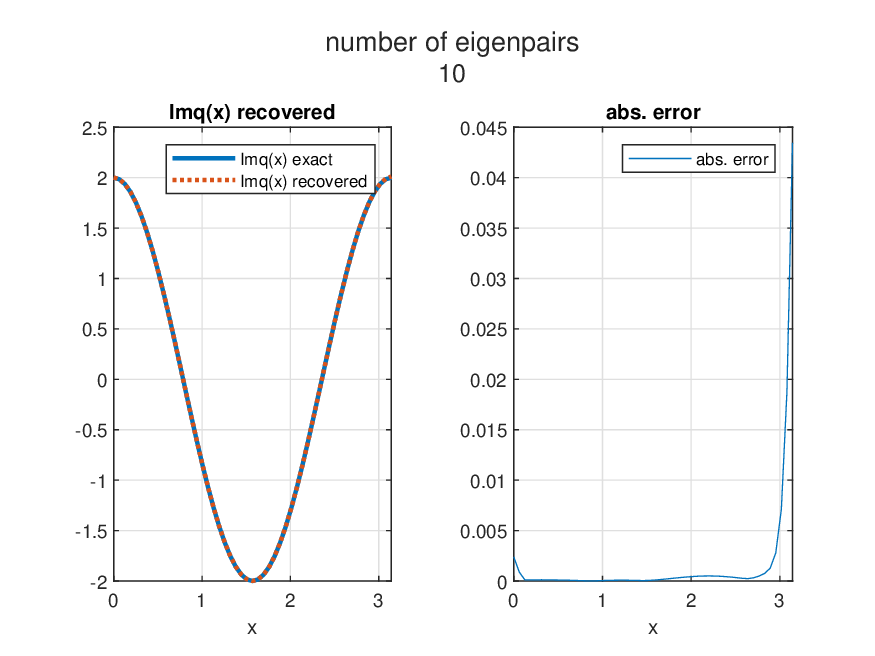}%
\caption{Mathieu potential from Example 4, recovered from ten noisy eigenpairs
with $\sigma=0.001$. The maximum absolute error resulted in $0.043$. The
constants $h=0.7$ and $H=i$ were recovered with the absolute error
$6\cdot10^{-5}$ and $3\cdot10^{-3}$, respectively.}%
\label{Fig9}%
\end{center}
\end{figure}

A greater level of noise: $\sigma=0.01$ still allows to recover the
Sturm-Liouville problem quite accurately, see Fig. 10. The maximum absolute
error resulted in $0.46$. The constants $h$ and $H$ were recovered with the
absolute error $3.12\cdot10^{-4}$ and $0.032$, respectively.%

\begin{figure}
[ptb]
\begin{center}
\includegraphics[
height=5.2112in,
width=6.9446in
]%
{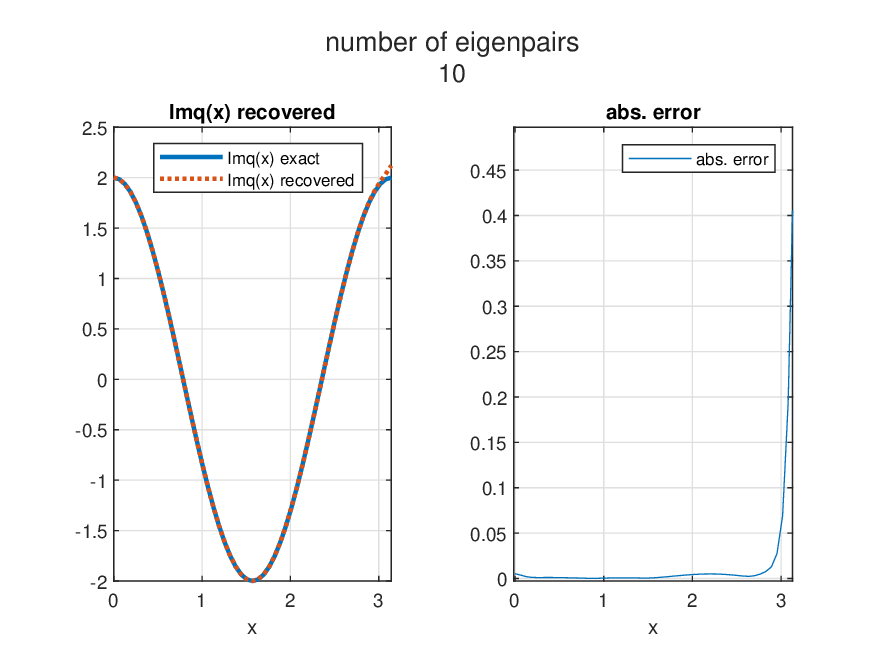}%
\caption{Mathieu potential from Example 4, recovered from ten noisy eigenpairs
with $\sigma=0.01$. The maximum absolute error resulted in $0.46$. The
constants $h$ and $H$ were recovered with the absolute error $3.12\cdot
10^{-4}$ and $0.032$, respectively.}%
\label{Fig10}%
\end{center}
\end{figure}
Even the level of noise $\sigma=0.1$ still allows one to obtain a meaningful
approximation, see Fig.11. The constants $h$ and $H$ were recovered with the
absolute error $0.0036$ and $0.32$, respectively.%

\begin{figure}
[ptb]
\begin{center}
\includegraphics[
height=5.1492in,
width=6.8616in
]%
{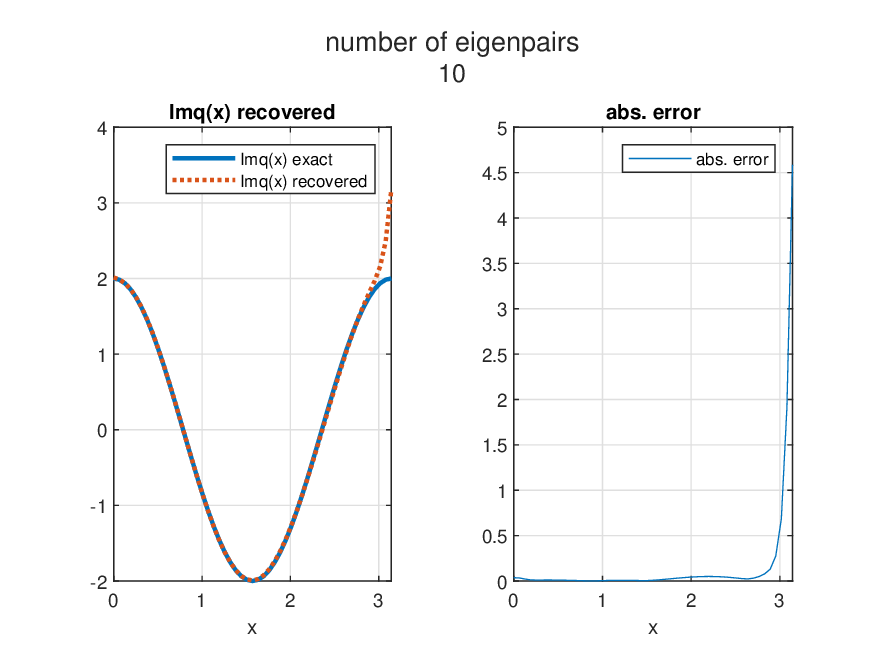}%
\caption{Mathieu potential from Example 4, recovered from ten noisy eigenpairs
with $\sigma=0.1$. The constants $h$ and $H$ were recovered with the absolute
error $0.0036$ and $0.32$, respectively.}%
\label{Fig11}%
\end{center}
\end{figure}

\subsection{Reconstruction from Weyl function, problem (IP2)}

The Weyl functions of the Sturm-Liouville problems considered in this
subsection were computed by the formula
\[
M(\rho)=-\frac{\Delta^{0}\left(  \rho\right)  }{\Delta\left(  \rho\right)  },
\]
where in their turn, $\Delta^{0}\left(  \rho\right)  $ and $\Delta\left(
\rho\right)  $ were computed with the aid of the NSBF representations
(\ref{phi NSBF}), (\ref{S NSBF}), (\ref{phi prime NSBF}), (\ref{S prime NSBF}%
), with the coefficients calculated following the recurrent integration
procedure from \cite{KNT} (see also \cite{KrBook2020}).

\textbf{Example 5. }Consider problem $L$ with the complex-valued potential and
complex constants in the boundary conditions:
\[
q(x)=e^{x}+\frac{i}{\left(  x+0.1\right)  ^{2}},\quad x\in\left[
0,\pi\right]  ,
\]%
\[
h=1-i,\quad H=e^{i}.
\]
Thus, the real and imaginary parts of the potential are the two Paine
potentials (see, e.g., \cite{Paine}, \cite{PryceBook}).

As the data for the inverse problem, the values of $M(\rho)$ were given at
$2000$ points logarithmically spaced along the segment $[0.01,1000]$, and
additionally at $20$ points $\rho_{j}$ distributed uniformy on the segment
$[0.01,1000]$. These last $20$ points were used for choosing an appropriate
value of the parameter $N_{1}$ in (\ref{sys Weyl}). Namely, $N_{1}$ was chosen
by minimizing the functional $Q(N_{1})=\max_{\rho_{j}}\left\vert \Delta
_{N_{1}}^{0}\left(  \rho\right)  +M(\rho)\Delta_{N_{1}}\left(  \rho\right)
\right\vert $.

The Weyl function computed on the segment $[0.01,100]$ is depicted in Fig. 12.%

\begin{figure}
[ptb]
\begin{center}
\includegraphics[
height=4.6218in,
width=6.1518in
]%
{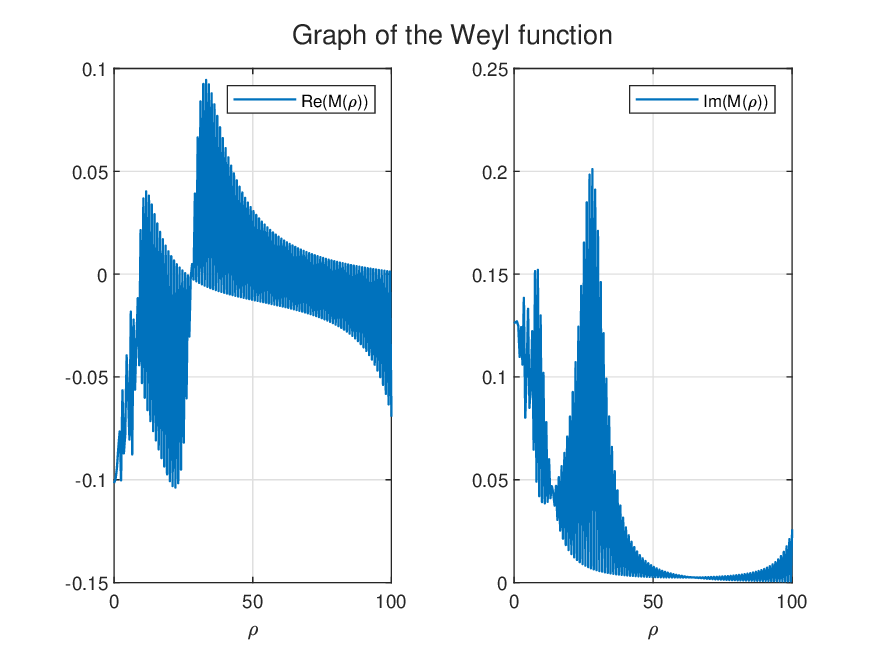}%
\caption{Weyl function of Sturm-Liouville problem from Example 5 on the
segment $[0.01,100]$.}%
\label{Fig12}%
\end{center}
\end{figure}

Fig. 13 presents the result of the reconstruction of the potential $q(x)$. The
maximum absolute error was $1.57$ attained at the origin (the relative error
was approximately $0.0157$). The constants $h$ and $H$ were recovered with the
absolute error $0.0044$ and $0.0013$, respectively. The parameter $N_{1}$ in
(\ref{sys Weyl}) resulted in $33$.%

\begin{figure}
[ptb]
\begin{center}
\includegraphics[
height=5.0214in,
width=6.6819in
]%
{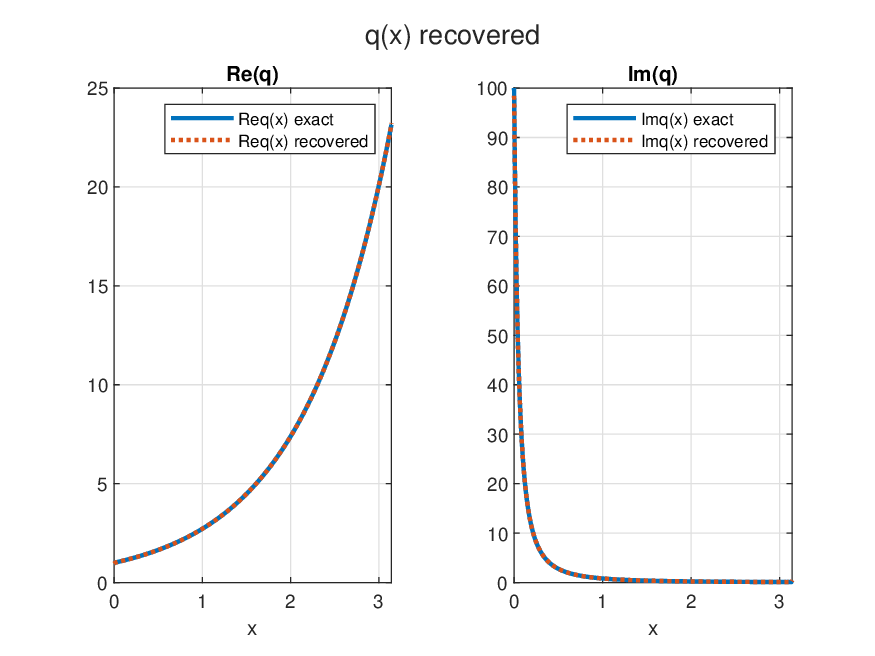}%
\caption{Potential from Example 5 recovered from $M(\rho)$ given at $2000$
points logarithmically spaced along the segment $[0.01,1000]$. The maximum
absolute error was $1.57$ attained at the origin (the relative error was
approximately $0.0157$). The constants $h$ and $H$ were recovered with the
absolute error $0.0044$ and $0.0013$, respectively. The parameter $N_{1}$ in
(\ref{sys Weyl}) resulted in $33$.}%
\label{Fig13}%
\end{center}
\end{figure}

\bigskip

\bigskip

\textbf{Example 6. }Consider problem $L$ with a less smooth potential
\[
q(x)=\left\vert 3-\left\vert x^{2}-3\right\vert \right\vert +i\left\vert
\cos\left(  2x\right)  \right\vert ,\quad x\in\left[  0,\pi\right]  ,
\]%
\[
h=e^{2i},\quad H=\pi-i.
\]
The Weyl function $M\left(  \rho\right)  $ was given at $2000$ points
logarithmically spaced on the segment $[0.01,1000]$.

Fig. 14 presents the result of the reconstruction of the potential $q(x)$. The
maximum absolute error resulted in $0.036$. The constants $h$ and $H$ were
recovered with the absolute error $7.4\cdot10^{-6}$ and $1.4\cdot10^{-4}$,
respectively. The parameter $N_{1}$ in (\ref{sys Weyl}) resulted in $49$. This
value was obtained by minimizing the value of the functional $Q(N_{1}%
)=\left\vert \Delta_{N_{1}}^{0}\left(  \rho\right)  +M(\rho)\Delta_{N_{1}%
}\left(  \rho\right)  \right\vert $ evaluated at 20 points $\rho_{j}$
distributed uniformy on the segment $[0.01,1000]$.%

\begin{figure}
[ptb]
\begin{center}
\includegraphics[
height=5.0807in,
width=6.7631in
]%
{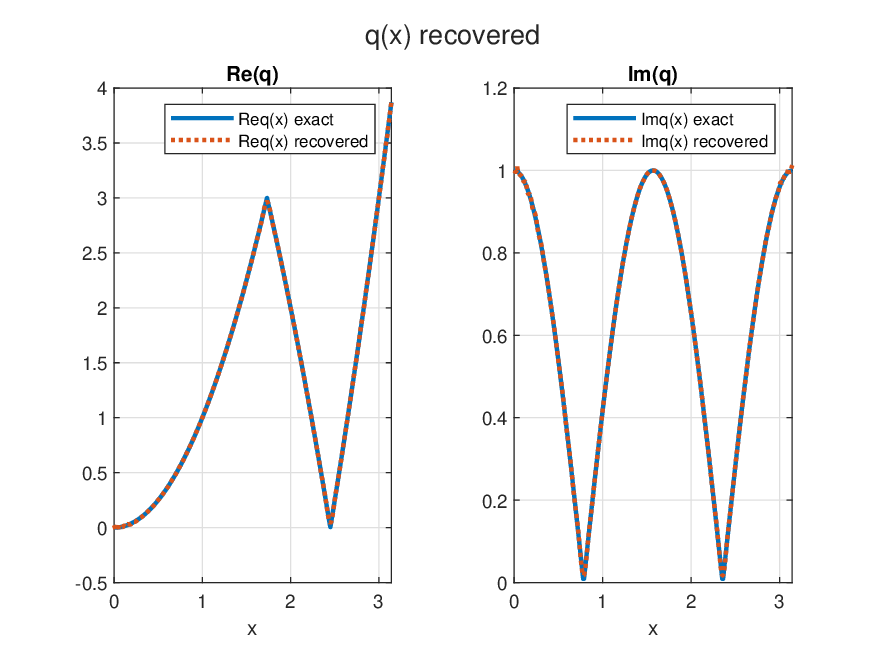}%
\caption{Potential from Example 6 recovered from the Weyl function $M\left(
\rho\right)  $ given at $2000$ points logarithmically spaced on the segment
$[0.01,1000]$. The maximum absolute error resulted in $0.036$. The constants
$h$ and $H$ were recovered with the absolute error $7.4\cdot10^{-6}$ and
$1.4\cdot10^{-4}$, respectively. }%
\label{Fig14}%
\end{center}
\end{figure}

\subsection{Reconstruction from eigenvalues and multiplier constants, problem
(IP3) \label{subsect Numeric IP3}}

The multiplier constants $\beta_{k}$ required as a part of the input data of
problem (IP3) are computed from (\ref{beta_k=}), where $\varphi_{h}(\rho
_{k},b)$ are computed by calculating the NSBF representations with the aid of
the recurrent integration procedure, as explained in \cite{KNT}.

\textbf{Example 7. }Consider the inverse problem (IP3) for the complex-valued
potential $q(x)$ and constants $h$ and $H$ from Example 3. Fig. 15 presents
the result of the recovery of $q(x)$ from 50 pairs $(\lambda_{k},\,\beta_{k}%
)$. The maximum absolute error resulted in $0.03$. The constants $h$ and $H$
were recovered with the absolute error $8.5\cdot10^{-5}$ and $6.21\cdot
10^{-5}$, respectively.%

\begin{figure}
[ptb]
\begin{center}
\includegraphics[
height=5.38in,
width=7.1691in
]%
{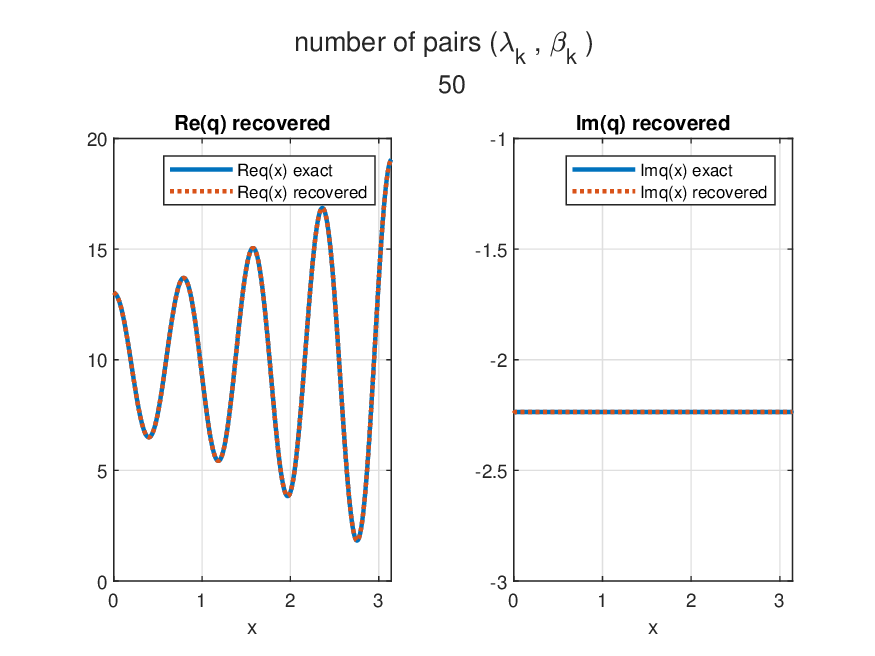}%
\caption{Complex valued potential from Example 3 recovered from 50 pairs
$(\lambda_{k},\,\beta_{k})$. The constants $h=\sqrt{2}$ and $H=-e$ were
recovered with the absolute error $8.5\cdot10^{-5}$ and $6.21\cdot10^{-5}$,
respectively.}%
\label{Fig15}%
\end{center}
\end{figure}

\textbf{Example 8. }Consider the inverse problem (IP3) for the complex-valued
potential $q(x)$ and constants $h$ and $H$ from Example 4. Fig. 16 presents
the result of the recovery of $q(x)$ from 10 pairs $(\lambda_{k},\,\beta_{k}%
)$. The maximum absolute error resulted in $0.0034$. The constants $h$ and $H$
were recovered with the absolute error $1.56\cdot10^{-4}$ and $1.2\cdot
10^{-4}$, respectively. The result is similar to that for problem (IP1),
compare to Fig. 8.%

\begin{figure}
[ptb]
\begin{center}
\includegraphics[
height=5.2158in,
width=6.951in
]%
{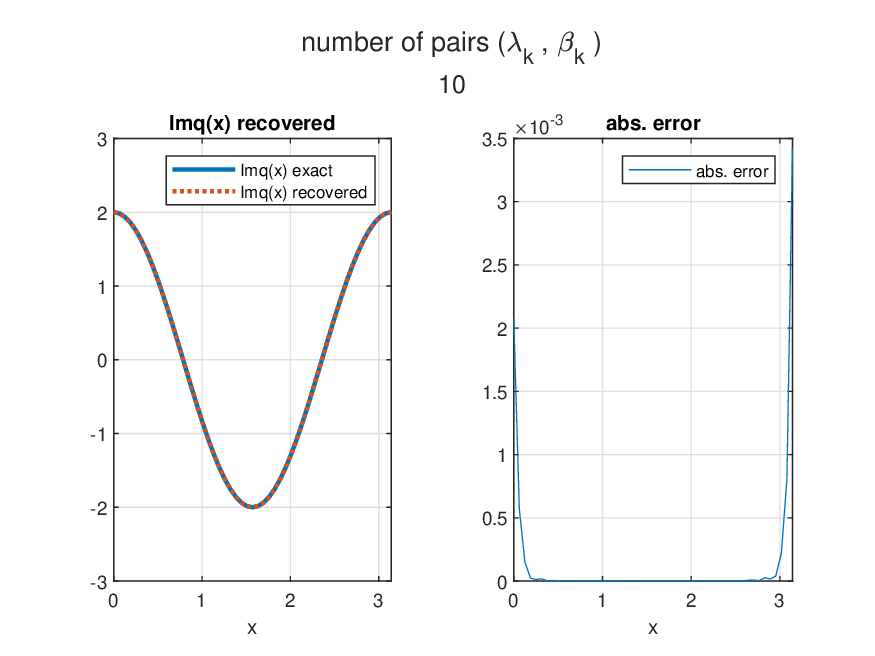}%
\caption{Complex valued Mathieu potential from Example 4 recovered from $10$
pairs $(\lambda_{k},\,\beta_{k})$. The maximum absolute error resulted in
$0.0034$. The constants $h$ and $H$ were recovered with the absolute error
$1.56\cdot10^{-4}$ and $1.2\cdot10^{-4}$, respectively. }%
\label{Fig16}%
\end{center}
\end{figure}

\subsection{Reconstruction from eigenvalues and norming constants, problem
(IP4)}

Since problem (IP4) reduces to problem (IP3) by computing first $\overset
{\cdot}{\Delta}(\rho_{k})$ and next the multiplier constants $\beta_{k}$, as
explained in subsection \ref{Subsect IP4}, here the main question is how
accurate is this reduction procedure. Let us consider an example. For the
Sturm-Liouville problem from Example 2 we computed the norming constants with
the aid of Matslise. This was done as follows. Matslise allows one to compute
normalized eigenfunctions of selfadjoint Sturm-Liouville problems $u_{k}(x)$:
$\left\Vert u_{k}\right\Vert _{%
\mathcal{L}%
^{2}\left(  0,b\right)  }=1$. Clearly, $\varphi_{h}(\rho_{k},x)=u_{k}%
(x)/u_{k}(0)$, and hence
\[
\alpha_{k}=\int_{0}^{b}\varphi_{h}^{2}(\rho_{k},x)dx=\frac{1}{u_{k}^{2}%
(0)}\int_{0}^{b}u_{k}^{2}(x)dx=\frac{1}{u_{k}^{2}(0)}.
\]
The values $u_{k}(0)$ \ were computed by Matslise.

We computed first fifteen norming constants for the Sturm-Liouville problem
from Example 2, and then with their aid, the values $\overset{\cdot}{\Delta
}(\rho_{k})$ and the multiplier constants $\beta_{k}$. We compared them to
those computed as explained in subsection \ref{subsect Numeric IP3}. The
maximum resulting difference was $2.8\cdot10^{-10}$. Thus, indeed, the
proposed procedure provides an accurate reduction of problem (IP4) to problem (IP3).

\section{Conclusions\label{Sect Conclusions}}

In this work a method for solving a variety of inverse spectral problems for
Sturm-Liouville equations with complex coefficients is presented. It is based
on the Neumann series of Bessel functions representations for solutions of the
Sturm-Liouville equations. The series representations possess several
remarkable features which make them especially convenient for solving inverse
problems. Most important, with respect to the spectral parameter $\rho
=\sqrt{\lambda}$ the series converge uniformly in any strip parallel to the
real axis, and the very first coefficients of the series contain all the
information on the Sturm-Liouville problem. Thus, the solution of the inverse
problems reduces to the solution of one or two systems of linear algebraic
coefficients for the coefficients of the series. The method is easy for
implementation, direct, accurate and applicable to a large variety of inverse
problems. Its performance is illustrated by numerical examples.\bigskip

\textbf{Acknowledgments} The author thanks Lady Estefania Murcia Lozano for
providing the eigenvalues for Example 4.

\textbf{Funding information }CONAHCYT, Mexico, grant \textquotedblleft Ciencia
de Frontera\textquotedblright\ FORDECYT - PRONACES/ 61517/ 2020.

\textbf{Data availability} The data that support the findings of this study
are available upon reasonable request.

\textbf{Conflict of interest }This work does not have any conflict of interest.

\end{document}